\documentclass[12pt]{article}
\usepackage{amsmath}
\usepackage{amssymb}
\usepackage{graphicx}
\usepackage[cp1251]{inputenc}
\usepackage[russian]{babel}
\newcommand{\il}[2]{\int\limits_{#1}^{#2}}
\newcommand{\ilp}[1]{\int\limits_{#1}^{+\infty}}

\newcommand{\ph}{\phantom{a}}
\newcommand{\phh}{\phantom{aaa}}
\newcommand{\ilpp}{\ilp{t_0}}
\newcommand{\sist}[2]{\left\{
\begin{array}{l}
{#1}\\
\ph\\
{#2}
\end{array}
\right.}

\newcommand{\und}[2]{\hskip -20pt
\begin{array}{l}
{_{^{#1}}}\\
{^{^{#2}}}
\end{array}
}

\begin{document}

MSC 39B99, 34K06, 34K11

\vskip 40pt

\centerline{\bf Oscillatory  criteria for the second order  linear}

 \centerline{\bf functional - differential equations with locally integrable coefficients}

\hskip 20 pt

\centerline{\bf G. A. Grigorian}

\hskip 20 pt

Abstract. The Riccati equation method is used to establish some oscillatory criteria for the second order linear functional - differential equation of multiple terms with locally integrable coefficients. An interval oscillation criterion for the second order linear functional - differential equation is proved. We have obtained a generalization of an oscillation theorem of L. Berezanski and E. Braverman, a generalization of the well known Fite's oscillation criterion and a  new global solvability criterion for the second order linear functional - differential equations with advanced and retarded arguments.

\hskip 20 pt

\vskip 20 pt

\centerline{\bf \S 1. Introduction}

\vskip 20 pt

Let  $q_k(t), \ph k=\overline{1,n}$  be locally integrable and  $\alpha_k(t), \ph k=\overline{1,n}$   be locally bounded real valued functions on  $[t_0;+\infty)$ and let  $p(t)$ be a positive function on  $[t_0;+\infty)$ such that $\frac{1}{p(t)}$ is locally integrable on  $[t_0;+\infty)$. Consider the equation
$$
(p(t)\phi'(t))' + \sum\limits_{k=1}^n q_k(t) \phi(\alpha_k(t)) = 0, \phh t\ge t_0. \eqno (1.1)
$$
Denote  $T_0 \equiv \min \{t_0; \hskip 2pt \inf_{\und{t\ge t_0}{k=\overline{1,n}}} \hskip -10 pt\alpha_k(t)\}$. A real valued continuous function  $\phi(t)$ on  $[T_0;+\infty)$   is called a solution of Eq. (1.1),  if $\phi'(t)$  is absolutely continuous on  $[t_0;+\infty)$  and $\phi(t)$ satisfies Eq. (1.1) almost everywhere on $[t_0;+\infty)$.

{\bf Definition 1.1.} {\it Eq. (1.1) is called oscillatory if its each solution has arbitrarily large zeroes.}

{\bf Definition 1.2} {\it Eq. (1.1) is called oscillatory on the segment  $[t_1;t_2] \linebreak(t_0\le t_1 < t_2 < + \infty)$, if its each solution vanishes on  $[t_1;t_2]$.}

The study of the question of oscillation  of linear functional - differential equations in particular of Eq. (1.1) is an important problem of qualitative theory of functional - differential equations and many works are devoted to it (see [1 - 6] and cited works therein). Among these works  note [1] where the following important result is proved.

{\bf Theorem 1.1}. {\it Suppose  $p(t) \equiv 1, \ph q_k(t) \ge 0, \ph 0 \le \alpha_k(t) \le t, \ph t\ge t_0, \linebreak \lim\limits_{t\to + \infty} \alpha_k(t) = + \infty, \ph k=\overline{1,n},$                                 and  for each $c =const > 0$  the ordinary differential equation
$$
\phi''(t) + \biggl(\sum\limits_{k=1}^n q_k(t) \frac{\alpha_k(t) - c}{t - c}\biggr) \phi(t) = 0, \phh t\ge t_0,
$$
is oscillatory. Then Eq. (1.1) is also  oscillatory.}

In this paper the Riccati equation method is used to establish oscillatory criteria for Eq. (1.1) in two new directions. The first direction is to obtain oscillatory criteria for Eq. (1.1) with both advanced and retarded arguments. The second one is to break nonnegativity condition imposed on the coefficients of Eq. (1.1) (see below Theorem 2.6 and Remark~ 2.1).
In the first direction we have obtained a generalization of Theorem 1.1  (see below  Theorem~2.4).
The first integral criterion of oscillation  for the equation
$$
\phi''(t) + q(t) \phi(t) = 0, \phh t\ge t_0, \eqno (1.2)
$$
where $q(t)$ is a continuous function on  $[t_0;+\infty)$
was formulated and proved by Fite (see [7]) stating  that if   $q(t)\ge 0, \ph t\ge t_0$, and  $\ilpp q(\tau) d\tau = + \infty,$  then Eq. (1.2) is oscillatory.
A generalization of this result is  Theorem 2.4  proved below.
In this work  an oscillation criterion on the finite segment for Eq. (1.1)is obtained (in the second direction; see below Theorem 2.6).

For studying of oscillatory property of Eq. (1.1) in this work  in general we shall letting the satisfaction of the following conditions (an exception is Theorem 2.6)

\noindent
A) \ph $q_k(t) \ge 0, \ph t\ge t_0, \ph k=\overline{1,n};$

\noindent
B) \ph $\lim\limits_{t\to +\infty} \alpha_k(t) = +\infty, \ph k=\overline{1,n}.$

The criteria of this work can be attributed to two groups.
To the first group we attribute the criteria with restriction

\noindent
C) \ph $\ilpp \frac{d\tau}{p(\tau)} = +\infty.$

To the second group we attribute the criteria without of restriction C).

\pagebreak

\centerline{\bf \S 2.  Main results}

\vskip 20 pt

Let $a_k(t) \phantom{a} k=\overline{1,n},$ be real valued locally integrable functions on $(-\infty;+\infty)$ and let $h_k \phantom{a} k=\overline{1,n},$ be some real constants. Consider the functional - differential equation
$$
\phi''(t) + \sum\limits_{k=1}^n a_k(t)\phi(t+h_k) = 0, \phantom{aaa} t\in(-\infty;+\infty) \eqno(2.1)
$$

A real valued continuous function $\phi(t)$ on $(-\infty;+\infty)$ is called a solution of Eq. (2.1) if $\phi'(t)$ is absolutely continuous on $(-\infty;+\infty)$ and $\phi(t)$ satisfies  Eq. (2.1) almost everywhere  on $(-\infty;+\infty)$. The restriction of any solution of Eq (2.1) on $[t_0;+\infty)$ we also will call a solution of Eq. (2.1).

{\bf Theorem 2.1.} {\it Suppose $|a_k(t)| \le a_k^0, \phantom{a} k=\overline{1,n}, \phantom{a} t\in(-\infty;+\infty)$, and let the transcen-\linebreak dent  equation
$$
\lambda^2 = \sum\limits_{k=1}^n a_k^0 e^{h_k\lambda} \eqno (2.2)
$$
has negative and positive solutions. Then for each $\gamma_0, \phantom{a} t_0 \in (-\infty; +\infty)$ Eq. (2.1) has the solution $\phi_0(t)$ on $(-\infty; +\infty)$ satisfying the conditions  $\phi_0(t_0) = \gamma_0, \phantom{a} \phi'_0(t_0) = 0$.}

Example 2.1. The equation
$$
\phi''(t) + \frac{1+ \varepsilon}{4(1 + t^2)}[\phi(t - 7/8) + \phi(t+2/3)] = 0 \eqno (2.3)
$$
for each $\gamma_0, \phantom{a} t_0 \in (-\infty; +\infty)$ and for enough small $\varepsilon > 0$  has the solution $\phi_0(t)$ on $(-\infty; +\infty)$ satisfying the conditions  $\phi_0(t_0) = \gamma_0, \phantom{a} \phi'_0(t_0) = 0$.

Example 2.2. The equation
$$
\phi''(t) + \frac{e \sin^2 t}{1 + e^2}\phi(t - 1) + \frac{1}{4(1 + t^2)}  \phi(t+1) = 0 \eqno (2.4)
$$
for each $\gamma_0, \phantom{a} t_0 \in (-\infty; +\infty)$  has the solution $\phi_0(t)$ on $(-\infty; +\infty)$ satisfying the conditions  $\phi_0(t_0) = \gamma_0, \phantom{a} \phi'_0(t_0) = 0$.

Example 2.3. The equation
$$
\phi''(t) + \frac{1}{7(1 + t^2)}[\phi(t - 3/2) +   \phi(t+3/2)]  = 0 \eqno (2.5)
$$
for each $\gamma_0, \phantom{a} t_0 \in (-\infty; +\infty)$  has the solution $\phi_0(t)$ on $(-\infty; +\infty)$ satisfying the conditions  $\phi_0(t_0) = \gamma_0, \phantom{a} \phi'_0(t_0) = 0$.

Example 2.4. The equation
$$
\phi''(t) + 32\sin t[\phi(t - 1/16) +   \phi(t+1/16)]  = 0 \eqno (2.6)
$$
for each $\gamma_0, \phantom{a} t_0 \in (-\infty; +\infty)$  has the solution $\phi_0(t)$ on $(-\infty; +\infty)$ satisfying the conditions  $\phi_0(t_0) = \gamma_0, \phantom{a} \phi'_0(t_0) = 0$.

Denote: $\Omega \equiv \{1, \dots, n\}$ \ph ($n\ge 1$), \ph $\Omega_+\equiv\{k\in\Omega: \alpha_k(t)\ge~t, \ph t\ge t_k,$  for some $t_k \ge t_0\}, \ph \Omega_-\equiv\{k\in \Omega: \alpha_k(t) \le t, \ph t\ge t_k,$ for some  $t_k \ge t_0\}$. Consider the equations
$$
(p(t)\phi'(t))' + \Bigl[\sum\limits_{k\in \Omega_\pm} q_k(t)\
\Bigr] \phi(t) =0, \ph t\ge t_0. \eqno (2.7_\pm)
$$

{\bf Theorem 2.2.} {\it  Let the conditions A)  and  B) be satisfied, and let $\Omega_\pm \ne \emptyset$ the equations  $(2.7_\pm)$  be oscillatory.  Then Eq. (1.1)  is oscillatory.}

In what follow under the symbol $\il{\xi}{\eta} u(\tau) d\tau$ we shall mean the integration of the function $u(t)$ by the direction from $\xi$ to $\eta$.

{\bf Theorem 2.3.} {\it Let the conditions A) and   B) be satisfied, $\Omega_\pm \ne \emptyset$   and let

\noindent
D) \ph $\ilpp\Bigl[\sum\limits_{k\in\Omega_+}q_k(\tau)\Bigr] d \tau < + \infty$.

\noindent
Then if the equations $(2.7_-)$  and
$$
(p(t)\phi'(t))' + \biggl[\sum\limits_{k\in\Omega_+}q_k(t)\exp\Bigl\{\il{t}{\alpha_k(t)} \frac{d \tau}{p(\tau)}\ilp{\tau}\Bigl(\sum\limits_{k\in\Omega_+}q_k(s)\Bigr) d s\Bigr\}\biggr] \phi(t) = 0, \phh t\ge t_0, \eqno (2.8)
$$
are oscillatory then Eq. (1.1) is also  oscillatory.}

The next theorem is a generalization of the mentioned above Fite's result ($p(t) \equiv~ 1, \linebreak n=~1, \phantom{a} \alpha_1(t) \equiv t$).

{\bf Theorem 2.4.} {\it Let the conditions  A)  - C) be satisfied  and let

\noindent
Е) \ph $\ilpp\Bigl(\sum\limits_{k=1}^n q_k(\tau)\Bigr) d \tau = +\infty.$

\noindent
Then Eq. (1.1) is oscillatory.}

Consider the equation
$$
(p(t)\phi'(t))' + \Biggl[\sum\limits_{k\in\Omega_+}q_k(t)\exp\Bigl\{\il{t}{\alpha_k(t)} \frac{d \tau}{p(\tau)}\ilp{\tau}\Bigl(\sum\limits_{k\in\Omega_+}q_k(s)\Bigr) d s\Bigr\} + \phantom{aaaaaaaaaaaaaaaaaaaaaaaaaaaa}
$$
$$
\phantom{aaaaaaaaaaaaaaaa}+\sum\limits_{k\in \Omega_- \backslash \Omega_+} q_k(t) \frac{\il{t_1}{\alpha_k(t)}\frac{d \tau}{p(\tau)}}{\il{t_1}{t}\frac{d \tau}{p(\tau)}}\Biggr] \phi(t) = 0, \phh t > t_1 \ge t_0, \eqno (2.9)
$$
The next result is a generalization of Theorem 1.1 ($p(t) \equiv 1, \phantom{a} \Omega_+ = \emptyset$)

{\bf Theorem 2.5.} {\it  Let the conditions A) - D) be satisfied, $\Omega_ - \cup \Omega_+ \ne \emptyset$   and let for each  $t_1 \ge t_0$  Eq. (2.9) be oscillatory. Then Eq. (1.1) is oscillatory.}

Consider the equation
$$
(p(t)\phi'(t))' + \Biggl[\sum\limits_{k\in\Omega_+}q_k(t)\exp\Bigl\{\il{t}{\alpha_k(t)} \frac{d \tau}{p(\tau)}\ilp{\tau}\Bigl(\sum\limits_{k\in\Omega_+}q_k(s)\Bigr) d s\Bigr\} + \phantom{aaaaaaaaaaaaaaaaaaaaaaaaaaaa}
$$
$$
\phantom{aaaaaaaaaaaaaaaaaaaaaaaaaa}+\sum\limits_{k\in \Omega_- \backslash \Omega_+} l_k \hskip 2pt q_k(t) \frac{\il{t_0}{\alpha_k(t)}\frac{d \tau}{p(\tau)}}{\il{t_0}{t}\frac{d \tau}{p(\tau)}}\Biggr] \phi(t) = 0, \phh t > t_0. \eqno (2.10)
$$

{\bf Corollary 2.1} {\it Let the conditions  A)   -  D) be satisfied, $\Omega_ - \cup \Omega_+ \ne \emptyset$  and let for some  $l_k\in (0,1), \phantom{a} k\in \Omega_-\backslash \Omega_+$,  Eq. (2.10) be oscillatory. Then Eq. (1.1) is   oscillatory.}

Let  $t_0 \le t_1 < t_2 \le t_3 < t_4 < + \infty$. Denote: $\omega_+\equiv \{k\in\Omega : \alpha_k(t) \ge t, \linebreak t\in [t_1,t_2]\}; \ph \omega_1^-\equiv\{k\in \Omega : t_1\le \alpha_k(t) \le t, \ph t\in [t_1,t_2]\}; \ph \omega_2^-\equiv \{k\in \Omega : \alpha_k(t) \le t, \linebreak t\in[t_3, t_4]\}; \ph T_1\equiv \inf_{ \und{t\in [t_1,t_4]}{k\in \Omega}} \alpha_k(t); \ph T_2\equiv \sup_{ \und{t\in [t_1,t_4]}{k\in \Omega}} \alpha_k(t); \ph t_2^+ \equiv \sup_{ \und{t\in [t_1,t_2]}{k\in \omega_+}} \alpha_k(t); \linebreak t_3^- \equiv \inf_{ \und{t\in [t_3,t_4]}{k\in \omega_2^-}} \alpha_k(t).$
Consider the equations
$$
(p(t)\phi'(t))' + \Biggl[\sum\limits_{k\in\omega_+}q_k(t)\exp\Bigl\{\il{t}{\alpha_k(t)} \frac{d \tau}{p(\tau)}\il{\tau}{t_2}\Bigl(\sum\limits_{k\in\omega_+}q_k(s)\Bigr) d s\Bigr\} + \phantom{aaaaaaaaaaaaaaaaaaaaaaaaaaaa}
$$
$$
\phantom{aaaaaaaaaaaaaaaaaaa}+\sum\limits_{k\in \omega_1^-} q_k(t) \frac{\il{t_1}{\alpha_k(t)}\frac{d \tau}{p(\tau)} + \varepsilon}{\il{t_1}{t}\frac{d \tau}{p(\tau)}+ \varepsilon}\Biggr] \phi(t) = 0, \phh t \in [t_1, t_2], \phh \varepsilon > 0; \eqno (2.11)
$$
$$
(p(t)\phi'(t))' + \Bigl[\sum\limits_{k\in\omega_2^-} q_k(t)\Bigr] \phi(t) = 0, \phh t\in [t_3,t_4]. \eqno (2.12)
$$

{\bf Theorem 2.6.} {\it Let $\omega_+ \cup \omega_1^- \ne \emptyset, \phantom{a} \omega_2^- \ne \emptyset$ and let the following conditions be satisfied:

\noindent
a) \ph $q_k(t) \ge 0, \ph t\in [T_1, T_2], \ph k \in\Omega;$

\noindent
b) \ph $t_2^+ \le t_3^-$;

\noindent
c) \ph for some $\varepsilon_0 > 0$  Eq. (2.11) is oscillatory on $[t_1,t_2]$  for all  \ph $\varepsilon \in (0,\varepsilon_0)$;

\noindent
d) \ph Eq. (2.12) is oscillatory on  $[t_3, t_4]$.

\noindent
Then Eq. (1.1) is oscillatory on $[T_1,T_2]$.}

{\bf Remark 2.1}. {\it If for Eq. (1.1) the conditions of Theorem 2.6 are fulfilled on each of intervals $[L_{2m};L_{2m+1}], \phantom{a} m=1,2,\dots,$ where $L_1 < L_2 \dots L_m \dots$ and $\lim\limits_{m\to+\infty}L_m =+\infty$ then Eq. (1.1) is oscillatory. Note that in this case outside of the set $\cup_{m=1}^{+\infty}[L_{2m};L_{2m+1}]$ the functions  $q_k(t), \phantom{a} k = \overline{1,n}$, can change their signs.}

Example 2.4. By Theorem 2.2 Eq. (2.3) is oscillatory.

Example 2.5. By Theorem 2.4 Eq. (2.4) is oscillatory.

Example 2.6. By Theorem 2.5 Eq. (2.5) is oscillatory.

Example 2.7. Show that Eq. (2.6) is oscillatory. According to Remark 2.1 it is enough to show that Eq. (2.6) is oscillatory on each interval $[2m\pi;(2m +1)\pi, \phantom{a} m=1,2,\dots$. Let $m$ be fixed. set $t_1 = (2m + \frac{1}{6})\pi, \phantom{a} t_2 = (2m + \frac{1}{2})\pi - \frac{1}{16}, \phantom{a} t_3 = (2m + \frac{1}{2})\pi + \frac{1}{16}, t_4= (2m +1)\pi  - \frac{\pi}{6}$. Then we have $T_1 = (2m + \frac{1}{6})\pi - \frac{1}{16},\phantom{a} T_2 = (2m +\frac{5}{6})\pi + \frac{1}{16}, \phantom{a} t_2^+ = t_3^- = (2m + \frac{1}{2})\pi.$ From here it follows that for Eq. (2.6) the conditions a) and b) of Theorem 2.6 are fulfilled. By Sturm's comparison
$$
\phi''(t) + 32\sin t \phi(t) = 0
$$
ia oscillatory on each of intervals $[t_1;t_2]$ and $[t_3;t_4]$. Then (by Sturm's comparison theorem) the conditions c) and d) of Theorem 2.6 for Eq. (2.6) are fulfilled. Therefore Eq. (2.6) is oscillatory on $[T_1;T_2] (\subset [2m\pi;(2m+1)\pi])$.

Example 2.8. Consider the equation
$$
\phi''(t) + \frac{\lambda}{t \ln t \ln (\ln t)} \phi( \ln t) = 0, \phh t\ge 3, \ph \lambda > 0 \eqno (2.13)
$$
It is not difficult to verify that the equation
$$
\phi''(t) + \frac{\lambda (\ln t - t_1)}{t \ln t \ln (\ln t) (t - t_1)} \phi(t) = 0, \ph t \ge 3,
$$
is not oscillatory for all  $\lambda > 0, \ph t > t_1$ (this fact follows from the inequality $\frac{\lambda (\ln t - t_1)}{t \ln t \ln (\ln t) (t - t_1)} \le \frac{1}{4 t^2}$ for all enough large t and from the non oscillation of Euler's equation $\phi''(t) + \frac{1}{4 t^2}\phi(t) = 0$ via Sturm's comparison theorem). Therefore Theorem 2.5  is not applicable to Eq.  (2.13).  Since for each  $\lambda > 0 \ph \ilp{3}\frac{\lambda d t}{t \ln t \ln(\ln t)} = + \infty$,  by Theorem 2.4 Eq. (2.13) is oscillatory.

\vskip 20 pt

\centerline{ \bf \S 3. Proof of the main results}

\vskip 20 pt

Let  $a(t)$ and  $b(t)$  be real valued locally integrable functions on $[t_0;+\infty)$. Consider the Riccati equation
$$
y'(t) + a(t) y^2(t) + b(t) = 0, \phh t\ge t_0. \eqno (3.1)
$$
An absolutely continuous function  $y_0(t)$   on  $[t_1;t_2) \ph (t_0 \le t_1 < t_2 \le +\infty)$  is called a solution of Eq. (3.1) on  $[t_1;t_2)$, if  $y_0(t)$  satisfies
(3.1) almost everywhere on  $[t_1;t_2)$.

{\bf Definition 3.1.} {\it  The set  $[t_1;t_2)$ is called maximal existence interval for the solution   $y(t)$ of Eq. (3.1) on  $[t_1;t_2)$,  if  $y(t)$   cannot be continued to the right of  $t_2$ as a solution of Eq. (3.1).}

The solutions  $y(t)$ of Eq. (3.1), existing on  $[t_1;t_2)$, are connected with the solutions    $(\phi(t), \psi(t))$   of the system of equations
$$
\sist{\phi'(t) = \phh \phh \phh a(t) \psi(t);}{\psi'(t) = - b(t) \phi(t), \ph t\ge t_0.} \eqno (3.2)
$$
by equalities  (see  [8], pp. 153 - 154):
$$
\phi(t) = \phi(t_1) \exp\biggl\{\il{t_1}{t} a(\tau) y(\tau) d \tau\biggr\}, \ph \phi(t_1) \ne 0, \phh \psi(t) = y(t) \phi(t). \eqno (3.3)
$$
Under a solution $(\phi(t), \psi(t))$ of the system (3.2) on any interval $[t_1;t_2) (\subset [t_0;+\infty))$ we mean an ordered pair of absolutely continuous functions $\phi(t), \phantom{a} \psi(t)$, defined on $[t_1;t_2)$, satisfying (3.2) almost everywhere on  $[t_1;t_2)$. Using contraction mapping method it is not difficult to show that for each values $\alpha_0$ and $\beta_0$ the system (3.2) has the unique solution $(\phi(t), \psi(t))$ on $[t_0; +\infty)$, satisfying the initial conditions $\phi(t_0) = \alpha_0, \phantom{a} \psi(t_0) = \beta_0$.
Hereafter we will assume that all solutions of equations and systems of equations be real valued and all functional equalities and inequalities which are considering on some sets we will understand their fulfillment almost everywhere on these  sets.

Let  $b_1(t)$ be a  real valued locally integrable function on  $[t_0;+\infty)$. Along with Eq.  (3.1) consider the Riccati equation
$$
y'(t) + a(t) y^2(t) + b_1(t) = 0, \phh t\ge t_0. \eqno (3.4)
$$
Before starting to prove the main results we need to prove two lemmas

{\bf Lemma 3.1}. {\it  Let  $a(t) \ge 0, \ph b(t) \le b_1(t), \ph t\in [t_1;t_2) \ph (t_0 \le t_1 < t_2 \le + \infty)$, and let  $y_1(t)$ be a solution of Eq. (3.4) on  $[t_1;t_2)$. Then each solution  $y_0(t)$   of Eq. (3.1) with  $y_0(t_1) \ge y_1(t_1)$ exists on  $[t_1;t_2)$,  and
$$
y_0(t) \ge y_1(t), \ph t\in [t_1;t_2). \eqno (3.5)
$$
}

Proof. Let  $[t_1;T)$  be the maximal existence interval for  $y_0(t)$. Show that  $T \ge t_2$. Suppose  $T < t_2$.  Then by virtue of (3.1) and  (3.4) we have:
$$
(y_0(t) - y_1(t))' + a(t) [y_0(t) + y_1(t)] (y_0(t) - y_1(t)) + b(t) - b_1(t) = 0, \phh t\in [t_1; T).
$$
From here is seen that  $y_0(t) - y_1(t)$  is a solution to the linear equation
$$
Y'(t) + a(t) [y_0(t) + y_1(t)] Y(t) + b(t) - b_1(t) = 0, \phh t\in [t_1;T).
$$
Therefore according to Cauchy formula we have:
$$
y_0(t) - y_1(t) = \exp\Bigl\{- \il{t_1}{t} a(\tau)\bigl[y_0(\tau) + y_1(\tau)\bigr]d\tau\Bigr\}\Biggl[y_0(t_1) - y_1(t_1) +\phantom{aaaaaaaaaaaaaaaaaaaaaaa}
$$
$$
\phantom{aaaaaaaaaaaaaaaaaaa}+ \il{t_1}{t}\exp\Bigl\{\il{t_1}{\tau} a(s) \bigl[y_0(s) + y_1(s)\bigr]d s\Bigr\}\bigl(b_1(\tau) - b(\tau)\bigr)d \tau\Biggr], \phh t\in [t_1;T).
$$
From here and from the conditions of the lemma it follows that
$$
y_0(t) \ge y_1(t), \phh t\in [t_1;T). \eqno (3.6)
$$
Let  $\phi_0(t) \equiv \exp\bigl\{\il{t_1}{t} a(\tau) y_0(\tau) d\tau\bigr\}, \ph \psi_0(t) = y_0(t)\phi_0(t), \ph t\in [t_1;T).$  By  (3.1)  - (3.3) \linebreak$(\phi_0(t), \psi_0(t))$ is a solution of the system (3.2) on $[t_1;T)$, which can be continued on  $[t_0;+\infty)$ as a solution of Eq. (3.2).
From (3.6) and from the inequalities   $a(t) \ge 0, \linebreak t\in [t_1;t_2), \ph T < t_2$ it follows that  $\phi_0(t) \ne 0, \ph t\in [t_1;t_3)$,  for some  $t_3 > T$. Then by  (3.3) $\widetilde{y}_0(t) \equiv \frac{\psi_0(t)}{\phi_0(t)}$  is a solution of Eq. (3.1) on  $[t_1;t_3)$. Obviously   $\widetilde{y}_0(t)$    coincides with  $y_0(t)$   on  $[t_1;T)$.   Therefore  $[t_1;T)$ is not the maximal existence interval for   $y_0(t)$. The obtained contradiction shows that  $T \ge t_2.$
This inequality with (3.6) prove (3.5). The lemma is proved.

Let $\phi_0(t)$  be a solution of Eq. (1.1) such that  $\phi_0(t) \ne 0, \ph t\in [t_1;t_2) \ph (t_0 \le t_1 < t_2 \le + \infty)$. Then it is not difficult to check that for the function  $y_0(t) \equiv \frac{p(t)\phi_0'(t)}{\phi_0(t)}$  the following equality takes place
$$
y_0'(t) + \frac{1}{p(t)} y_0^2(t) + \sum\limits_{k=1}^n q_k(t) \exp\biggl\{\il{t}{\alpha_k(t)}\frac{y_0(\tau)}{p(\tau)} d\tau\biggr\} = 0, \phh t\in[\widetilde{t}_1;\widetilde{t}_2), \eqno (3.7)
$$
 where $[\widetilde{t}_1;\widetilde{t}_2)\ph (\subset [t_1;t_2))$  is defined from the condition: $\phi_0(\alpha_k(t)) \ne 0, \ph t\in [t_1;t_2), \linebreak k= \overline{1,n}$.

{\bf Remark 3.1}. {\it  By (1.1) the function  $y_0(t)$ is absolutely continuous on  $[t_1;t_2)$.

{\bf Lemma 3.2}. {\it Let the conditions  A) - C) be fulfilled, and let  $\phi_0(t)$  be a solution of Eq.  (1.1) such that  $\phi_0(t) \ne 0, \ph t\ge t_1$, for some  $t_1 \ge t_0$.
Then the function  $y_0(t) \equiv \frac{p(t)\phi_0'(t)}{\phi_0(t)}, \linebreak t\ge t_1,$ is nonnegative on $[t_2;+\infty)$, where $t_2\ge t_1$  such that $\alpha_k(t) \ge t_1$ for  $t\ge t_2, \ph k= \overline{1,n}.$}

Proof. By  B)  chose  $t_2 \ge t_1$  such that  $\alpha_k(t) \ge t_1$  for  $t\ge t_2, \ph k= \overline{1,n}.$  Then by virtue of (3.7) $y_0(t)$ is a solution of the Riccati Equation
$$
y'(t) + \frac{1}{p(t)} y^2(t) + \sum\limits_{k=1}^n q_k(t) \exp\biggl\{\il{t}{\alpha_k(t)}\frac{y_0(\tau)}{p(\tau)} d\tau\biggr\} = 0, \phh t \ge t_2.  \eqno (3.8)
$$
Suppose that for some $t_3 \ge t_2, \ph y_0(t_3) < 0$. Along with Eq. (3.8) consider the equation
$$
y'(t) + \frac{1}{p(t)} y^2(t) = 0, \phh t\ge t_3. \eqno (3.9)
$$
Let  $y_1(t)$ be the solution of this equation with  $y_1(t_3) = y_0(t_3)$. Then obviously
$$
y_1(t) = \frac{1}{\il{t_3}{t}\frac{d\tau}{p(\tau)} + \frac{1}{y_0(t_3)}}, \phh t\in [t_3;t_4),
$$
where  $[t_3;t_4)$ is the maximal existence interval for  $y_1(t)$. From the condition  C) it follows that  $ \il{t_3}{t_4}\frac{d\tau}{p(\tau)} + \frac{1}{y_0(t_3)} = 0$ and $t_4 < +\infty$. On the other hand
applying  Lemma 3.1  to the equations (3.8), (3.9)  and  (in view of  A)) taking into account the inequality \linebreak$\sum\limits_{k=1}^n q_k(t) \exp\biggl\{\il{t}{\alpha_k(t)} \frac{y_0(\tau)}{p(\tau)}d \tau\biggr\} \ge~0, \ph t \ge t_3$,  we conclude that  $y_1(t)$ exists on  $[t_3;+\infty)$.  The obtained contradiction completes the proof of the lemma.

{\bf 3.1. Proof of Theorem 2.1.} Let $\lambda_1 (< 0)$ and $\lambda_2 (> 0)$
the solutions od transcendent equation (y) Then for each real constants $c_1$ and $c_2$ the function $\chi_0(t) \equiv c_1 e^{\lambda_1 t} + c_2 e^{\lambda_2 t}$ is a solution of the equation
$$
\phi''(t) = \sum\limits_{k=1}^n a_k^0\phi(t + h_k) \eqno (3.10)
$$
on $(-\infty;+\infty)$. Chose $c_1 \ge 0$ and $c_2 \ge 0$ such that $c_1 + c_2 = |\gamma_0|, \phantom{a} c_1\lambda_1 + c_2\lambda_2 = 0$. Then we have $\chi_0(t_0) = |\gamma_0|, \phantom{a} \chi'_0(t_0) = 0$. This means that $\chi_0(t)$ is a solution of the functional - integral equation
$$
\chi(t) = |\gamma_0| + (A_0 \chi)(t), \phantom{aaa} t\in (-\infty;+\infty), \eqno (3.11)
$$
wherw $(A_0 \chi)(t)\equiv \int\limits_{t_0}^td\tau\int\limits_{t_0}^\tau\biggl(\sum\limits_{k=1}^n a_k^0\chi(s + h_k) d s\biggl)d s, \phantom{a} t\in (-\infty;+\infty)$. Therefore for $\chi_0(t)$ the equality
$$
\chi_0(t) = |\gamma_0| + \sum_{m=1}^N(A_0^mf_0)(t) + (A_0^{N+1}\chi_0)(t), \phantom{aaa} t\in (-\infty;+\infty),
$$
takes place for arbitrary $N=1,2, \dots$, where $f_0(t) \equiv |\gamma_0|, \phantom{a} t\in (-\infty;+\infty)$. Hence
$$
 0 \le |\gamma_0| + \sum_{m=1}^N(A_0^mf_0)(t) \le \chi_0(t), \phantom{aaa} t\in (-\infty;+\infty),
$$
From here it follows that the series $F_0(t) \equiv |\gamma_0| + \sum_{m=1}^{+\infty}(A_0^mf_0)(t) \phantom{a} t\in (-\infty;+\infty)$ converges for each $t\in (-\infty;+\infty)$. Denote: $(A f)(t)\equiv \int\limits_{t_0}^td\tau\int\limits_{t_0}^\tau\biggl(\sum\limits_{k=1}^n a_k^0 f(s + h_k) d s\biggl)d s, \phantom{a} t\in (-\infty;+\infty)$, for arbitrary continuous function $f(t)$ on  $(-\infty;+\infty)$. Taking into account the conditions of the theorem it is not difficult to verify that
$$
|(A^m f_0)(t)| \le (a_0^m f_0)(t), \phantom{aaa}m=1,2,\dots, \phantom{aaa}t\in (-\infty;+\infty).
$$
From here and from the convergence of $F_0(t)$ it follows the convergence of the series
$$
\phi_0(t) \equiv \gamma_0 + \sum_{m=1}^{+\infty}(A^mf_0)(t), \phantom{aaa} t\in (-\infty;+\infty).
$$
Show that $\phi_0(t)$ is the required solution. For this it is enough to sow that $\phi_0(t)$ is a solution of the functional - integral equation
$$
\phi(t) = \gamma_0 + (A\phi)(t), \phantom{aaa} t\in (-\infty;+\infty). \eqno(3.12)
$$
Denote: $F_N(t) \equiv \gamma_0 + \sum\limits_{m=1}^N(A^mf_0)(t), \phantom{a} t\in (-\infty;+\infty).$ We have $\gamma_0 +(AF_N)(t) - F_N(t) = (A^{N+1} f_0)(t) \to 0$ uniformly  on $[-T;T]$ for each $T>0$. Hence
$$
\phi_0(t) = \lim\limits_{N\to +\infty} F_N(t) = \gamma_0 + \lim\limits_{N\to +\infty}(AF_n)(t) \eqno (3.13)
$$
uniformly on $[-T;T]$. For arbitrary $\varepsilon >0$ we can chose $N=N(\varepsilon)$ so large that $|(A\phi_0)(t) - (F_N(t)| \le \varepsilon, \phantom{a} t\in [-T;T]$. It means that $\lim\limits_{N\to +\infty} (AF_n)(t) = (A\phi_0)(t)$ uniformly on $[-T;T]$. Since $T > 0$ is arbitrary from here and from (6) it follows that $\phi_0(t)$ is a solution to Eq. (5) on $(-\infty;+\infty)$. The theorem is proved.

{\bf 3.2. Proof of Theorem 2.2.}  Suppose some solution  $\phi_1(t)$ of Eq. (1.1) has no arbitrarily  large zeroes. Let then  $\phi_1(t) \ne 0, \ph t\ge t_1$  for some  $t_1 \ge t_0$. Due to the condition  B)  chose  $t_2 \ge t_1$  so large that $\alpha_k(t) \ge t_1$  for $t \ge t_2, \ph k=\overline{1,n}$.  Then according to (3.7) for the function $y_1(t) \equiv \frac{p(t) \phi_1'(t)}{\phi_1(t)}, \ph t\ge t_1$,  the equality
$$
y_1'(t) + \frac{1}{p(t)} y_1^2(t) + \sum\limits_{k=1}^n q_k(t) \exp\biggl\{\il{t}{\alpha_k(t)}\frac{y_1(\tau)}{p(\tau)} d\tau\biggr\} = 0, \phh t \ge t_2.  \eqno (3.14)
$$
 is fulfilled. From here and from  A) it follows that  $y_1(t)$ is a monotonically non increasing function on $[t_2;+\infty)$. Then the following cases are possible.

\noindent
1} \ph $y_1(t) \ge 0, \ph t\ge t_2$;

\noindent
2) \ph $y_1(t) < 0, \ph t\ge t_3,$ for some  $t_3 \ge t_2$.

\noindent
Suppose the case 1) takes place.  Let $t_4 \ge t_2$ be such that    $\alpha_k(t) \ge t, \ph t\ge t_4, \ph k \in \Omega_+$. Then from A) it follows that
$$
\sum\limits_{k=1}^n q_k(t) \exp\biggl\{\il{t}{\alpha_k(t)}\frac{y_1(\tau)}{p(\tau)} d\tau\biggr\} \ge \sum\limits_{k\in \Omega_+} q_k(t), \phh t\ge t_4. \eqno (3.15)
$$
Consider the Riccati equations
$$
y'(t) + \frac{1}{p(t)} y^2(t) + \sum\limits_{k=1}^n q_k(t) \exp\biggl\{\il{t}{\alpha_k(t)}\frac{y_1(\tau)}{p(\tau)} d\tau\biggr\} = 0, \phh t \ge t_4.  \eqno (3.16)
$$
$$
y'(t) + \frac{1}{p(t)} y^2(t) + \sum\limits_{k\in \Omega_+}^n q_k(t) = 0, \phh t\ge t_4.  \eqno (3.17)
$$
By (3.14) $y_1(t)$ is a solution of Eq. (3.16) on  $[t_4;+\infty)$.  Let  $y_+(t)$ be the solution of Eq.  (3.17) with  $y_+(t_4) = y_1(t_4)$. Then applying  Lemma 3.1 to the equations (3.16) and  (3.17)
and taking into account (3.11) we conclude that $y_+(t)$ exists on  $[t_4;+\infty)$.
Hence  $\phi_+(t) \equiv \exp\biggl\{\il{t_4}{t}\frac{y_+(\tau)}{p(\tau)}\biggr\}, \ph t\ge t_4,$             is a solution to Eq. $(2.7_+)$ on  $[t_4;+\infty)$, which can be continued on $[t_0;+\infty)$  as a solution to Eq. (1.1).
It is evident that $\phi_+(t) \ne 0, \ph t\ge t_4$. Therefore  $(2.7_+)$  is not oscillatory which contradicts the condition of the theorem. Then it remains to suppose that   2) takes place. Let  $t_5 \ge t_3$ be  so large that $\alpha_k(t) \le t, \ph t\ge t_5, \ph k\in \Omega_-$. Then from A) it follows:
$$
\sum\limits_{k=1}^n q_k(t) \exp\biggl\{\il{t}{\alpha_k(t)}\frac{y_1(\tau)}{p(\tau)} d\tau\biggr\} \ge \sum\limits_{k\in \Omega_-} q_k(t), \phh t\ge t_5. \eqno (3.18)
$$
Consider the Riccati equation
$$
y'(t) + \frac{1}{p(t)} y^2(t) + \sum\limits_{k\in \Omega_-}^n q_k(t) = 0, \phh t\ge t_5.  \eqno (3.19)
$$
Let $y_-(t)$ be the solution of this equation with  $y_-(t_5) = y_1(t_5)$. Then applying Lemma~3.1
to the equations (3.16) and (3.19) and taking into account (3.18) we conclude that  $y_-(t)$ exists on  $[t_5;+\infty)$.   Hence  $\phi_-(t) \equiv \exp\biggl\{\il{t_5}{t}\frac{y_-(\tau)}{p(\tau)}\biggr\}, \ph t\ge t_5,$  is a solution of Eq. $(2.7_-)$ on  $[t_5;+\infty)$, which can be continued on  $[t_0;+\infty)$     as a solution of Eq. $(2.7_-)$. Obviously  $\phi_-(t) \ne 0, \ph t\ge t_5$. Therefore Eq. $(2.7_-)$ is not oscillatory, which contradicts  the condition of the theorem. The obtained contradiction completes the proof of the theorem.

{\bf 3.3. Proof of Theorem 2.3.} Let $\phi_1(t)$ be a solution of Eq. (1.1) such that  $\phi_1(t) \ne~0, \linebreak t\ge t_1$, for some $t_1\ge t_0$. Then for the function  $y_1(t) \equiv \frac{p(t)\phi_1'(t)}{\phi_1(t)}, \ph t\ge t_1,$  the equality (3.14) is satisfied.  Obviously $y_1(t)$  is a solution of the linear equation
$$
Y'(t) + \frac{1}{p(t)}y_1(t)Y(t) + \sum\limits_{k=1}^n q_k(t) \exp\biggl\{\il{t}{\alpha_k(t)}\frac{y_1(\tau)}{p(\tau)}\biggr\} = 0, \phh t\ge t_2,
$$
where $t_2 \ge t_1$  is so large that $\alpha_k(t) \ge t_1$ for  $t \ge t_2, \ph k=\overline{1,n}$  (existence of $t_2$ is guaranteed by the condition B)).
Therefore by Cauchy formula the equality
$$
y_1(t) = \exp\biggl\{ - \il{\xi}{t}\frac{y_1(\tau)}{p(\tau)} d\tau\biggr\}\biggl[y_1(\xi) - \il{\xi}{t}\Bigl(\sum\limits_{k=1}^n q_k(\tau) \exp\biggl\{\il{\xi}{\alpha_k(\tau)}\frac{y_1(s)}{p(s)} d s\biggr\}\Bigr)d\tau\biggr], \eqno (3.20)
$$
$t\ge \xi \ge t_2.$
is fulfilled.  By (3.14) $y_1(t)$ is a monotonically non increasing function on $[t_2;+\infty)$. Therefore for  $y_1(t)$   only the cases  1)  and  2) are possible.  Let the case 1) takes place.  Then from the conditions  A) and D) and from (3.19) we have:
$$
y_1(\xi) \ge \ilp{\xi}\biggl(\sum\limits_{k\in \Omega_+}q_k(\tau)\biggr) d \tau, \phh \xi \ge t_2.
$$
Hence
$$
\sum\limits_{k=1}^n q_k(t)\exp\biggl\{\il{t}{\alpha_k(t)}\frac{y_1(\tau)}{p(\tau)} d\tau\biggr\} \ge \sum\limits_{k\in \Omega_+} q_k(t)\exp\biggl\{\il{t}{\alpha_k(t)}\frac{y_1(\tau)d\tau}{p(\tau)}\ilp{\xi}\biggl(\sum\limits_{j\in \Omega_+}q_j(s)\biggr) d s\biggr\},  \eqno (3.21)
$$
$t\ge t_2.$
Further the arguments of the proof are analogous of the arguments of the proof of Theorem 2.2 only with the difference that in place of (3.15) is used (3.21). The proof of the theorem is complete.

{\bf 3.4. Proof of Theorem 2.4}. Suppose for some solution $\phi_1(t)$ of Eq. (1.1) the inequality  $\phi_1(t) \ne 0, \ph t\ge t_1$, is fulfilled for some  $t_1 \ge t_0$.
Then for the function \linebreak$y_1(t) \equiv \frac{p(t) \phi_1'(t)}{\phi_1(t)}, \ph t\ge t_1,$ the equality     (3.19) holds.  By virtue of Lemma 3.2 from the conditions  A) - C) it follows that $y_1(t) \ge 0, \ph t\ge t_2.$
Then from (3.20) it follows:
$$
\ilp{t_2}\biggl(\sum\limits_{k=1}^n q_k(\tau) \exp\biggl\{\il{t_2}{\alpha_k(\tau)}\frac{y_1(s)}{p(s)} d s \biggr\}\biggr) d\tau \le y_1(t_2). \eqno (3.22)
$$
Due to the condition B)  chose  $t_6 \ge t_2$ so large that  $\alpha_k(t) \ge t_2, \ph t\ge t_6, \ph k =\overline{1,n}.$ Then taking into account the conditions  A)  and  E) we will get:
$$
\ilp{t_2}\biggl(\sum\limits_{k=1}^n q_k(\tau) \exp\biggl\{\il{t_2}{\alpha_k(\tau)}\frac{y_1(s)}{p(s)} d s \biggr\}\biggr) d\tau \ge \ilp{t_2}\biggl(\sum\limits_{k=1}^n q_k(\tau)\biggr) d \tau = +\infty,
$$
which contradicts  (3.22). The obtained contradiction completes the proof of the theorem.

{\bf 3.5 Proof of Theorem 2.5.} Suppose for some solution $\phi_1(t)$ of Eq.  (1.1) the inequality   $\phi_1(t) \ne 0, \ph t\ge t_1,$  takes place for some  $t_1 \ge t_0$.
Then for the function  $y_1(t) \equiv \frac{p(t)\phi_1'(t)}{\phi_1(t)}, \linebreak t\ge~ t_1,$   the equality    (3.20) holds. By virtue of Lemma 3.1 from the conditions  A) -  C) it follows that  $y_1(t) \ge 0, \ph t\ge t_2.$ Then taking into account D) from (3.20) we will get:
$$
y_1(t) \ge \ilp{t} \biggl(\sum\limits_{k\in \Omega_+} q_k(\tau)\biggr) d \tau \phh t\ge t_2. \eqno (3.23)
$$
Let $\xi > t_2$. Consider the Riccati equations
$$
y'(t) + \frac{1}{p(t)} y^2(t) + \sum\limits_{k=1}^n q_k(t) \exp\biggl\{\il{t}{\alpha_k(t)} \frac{y_1(\tau)}{p(\tau)} d\tau\biggr\} = 0, \phh t\ge \xi; \eqno (3.24)
$$
$$
y'(t) + \frac{1}{p(t)} y^2(t) = 0,  \phh t\ge \xi; \eqno (3.25)
$$
Obviously  $y_\xi(t) \equiv 1/ \il{t_2}{t}\frac{d\tau}{p(\tau)}, \ph t \ge  \xi > t_2,$ is a solution to the last equation. Chose $\xi > t_2$ so close to  $t_2$, that $y_\xi(\xi) > y_1(\xi)$.  Then since $y_1(t)$ is a solution of Eq. (3.24) on $[\xi;+\infty)$,  by virtue of Lemma 3.1 from  A) it follows that
$$
y_1(t) \le 1/ \il{t_2}{t}\frac{d\tau}{p(\tau)},\ph t \ge  \xi. \eqno (3.26)
$$
Without loss of generality we will take that  $t_2 > t_0$ so large that   $\alpha_k(t)  \ge t, \ph k\in \Omega_+, \linebreak \alpha_k(t) \le ~t, \ph k\in \Omega_-, \ph t\ge t_2.$  Then taking into account the condition A) from (3.23) and (3.26) we will get:
$$
\sum\limits_{k=1}^n q_k(t) \exp\biggl\{\il{t}{\alpha_k(t)} \frac{y_1(\tau)}{p(\tau)} d \tau\biggr\} \ge \sum\limits_{k\in \Omega_+} q_k(t) \exp\biggl\{\il{t}{\alpha_k(t)} \frac{y_1(\tau)}{p(\tau)} d \tau\biggr\}+\phantom{aaaaaaaaaaaaaaaaaaaaaaaaa}
$$
$$
+ \sum\limits_{k\in \Omega_- \backslash \Omega_+} q_k(t) \exp\biggl\{\il{t}{\alpha_k(t)} \frac{y_1(\tau)}{p(\tau)} d \tau\biggr\} \ge \sum\limits_{k\in \Omega_+} q_k(t) \exp\biggl\{\il{t}{\alpha_k(t)} \frac{d\tau}{p(\tau)}\ilp{\tau}\biggl(\sum\limits_{j\in \Omega_+} q_j(s)\biggr) d s\biggr\} +
$$
$$
+ \sum\limits_{k\in \Omega_- \backslash \Omega_+} q_k(t) \frac{ \il{t_2}{\alpha_k(t)}\frac{d s}{p(s)}}{\il{t_2}{t}\frac{d s}{p(s)}} \stackrel{def}{=} Q(t), \phh \xi > t_2. \eqno (3.27)
$$
Consider the Riccati equation
$$
y'(t) + \frac{1}{p(t)} y^2(t) + Q(t) = 0, \phh t \ge \xi.
$$
Let $y_2(t)$  be the solution of this equation with  $y_2(\xi) = y_1(\xi)$.  Then since  $y_1(t)$ is a solution of Eq. (3.24) on  $[\xi;+\infty)$, by virtue of Lemma 3.1 from (3.27) it follows that  $y_2(t)$ exists on
$[\xi;+\infty)$.  Hence, $\phi_2(t) \equiv \exp\biggl\{\il{\xi}{t} \frac{y_2(\tau)}{p(\tau)} d \tau\biggr\}, \ph t\ge \xi,$  is a solution of Eq. (2.9)  for  $t_1 = t_2$  on  $[\xi; +\infty)$,  which can be continued on $[t_0;+\infty)$ as a solution of Eq. (2.9).  Obviously $\phi_2(t) \ne 0, \ph t\ge \xi.$ Therefore for   $t_1 = t_2$   Eq. (2.9) is not oscillatory, which contradicts the condition of the theorem. The obtained contradiction completes the proof of the theorem.

{\bf 3.6 Proof of Corollary 2.1}.  Let  $0 < l_k < 1, \ph k\in \Omega_-\backslash \Omega_+$.
Then from the \linebreak condition  C)  it follows that for each  $t_1 \ge t_0$ there exists  $t_2 > t_1$  such that  \linebreak $\il{t_1}{\alpha_k(t)}\frac{d s}{p(s)}/ \il{t_1}{t} \frac{d s}{p(s)}  \ge l_k \il{t_0}{\alpha_k(t)}\frac{d s}{p(s)}/ \il{t_0}{t} \frac{d s}{p(s)}, \ph t \ge t_2, \ph k\in \Omega_- $.  Therefore by virtue of the Sturm's comparison theorem ( see  [9], p. 334)  from the oscillation of Eq.  (2.9) it follows the oscillation of Eq. (2.8) for all  $t_1 > t_0$. The proof of the corollary is complete.

{\bf Remark 3.1.} {\it  From the proofs of Theorems 2.3  and  2.4 is seen that in their formulations the condition  C) can be replaced by condition

\noindent
$C')$ Eq.  $(2.7_-)$ is oscillatory.}

{\bf 3.7 Proof of Theorem 2.6.}  Suppose some solution  $\phi_1(t)$
of Eq. (1.1) does not vanish on  $[T_1;T_2]$. Then by  (3.7) the function  $y_1(t) \equiv \frac{p(t)\phi_1'(t)}{\phi_1(t)}, \ph t\in [T_1;T_2],$  is a solution of the Riccati equation
$$
y'(t) + \frac{1}{p(t)} y^2(t) + \sum\limits_{k=1}^n \exp\biggl\{\il{t}{\alpha_k(t)} \frac{y_1(\tau)}{p(\tau)} d\tau \biggr\} = 0, \phh t\ge t_0, \eqno (3.28)
$$
on  $[t_1;t_4]$. From here and from the condition  a) it follows that $y_1(t)$ is a monotonically non increasing  function on  $[t_1;t_4]$. Hence by virtue of the condition  b)  the cases

\noindent
$I) \ph y_1(t) \ge 0, \ph t\in [t_1; t_2^+], \phh II)\ph y_1(t) \le 0, \ph t\in [t_3^-;t_4],$

\noindent
are possible

\noindent
Let the case I) takes place. Consider the Riccatri equation
$$
y'(t) + \frac{1}{p(t)} y^2(t) = 0, \phh t\ge t_1. \eqno (3.29)
$$
It is evident that $y_\varepsilon(t) \equiv 1/ \Bigl(\il{t_1}{t} \frac{d\tau}{p(\tau)}+ \varepsilon\Bigr), \ph t\ge t_1, \ph \varepsilon > 0,$  is a solution of this equation on $[t_1;+\infty)$. It also is evident that for enough small  $\varepsilon_0 > 0$ the inequality
 $y_\varepsilon(t_1) \ge y_1(t_1)$ holds for all   $\varepsilon \in (0;\varepsilon_0)$.
Then applying Lemma 3.1 to the equations  (3.28) and  (3.29) and taking into account the condition  a) we conclude that
$$
0 \le y_1(t) \le \frac{1}{\il{t_1}{t}\frac{d\tau}{p(\tau)} + \varepsilon}, \phh t\in[t_1;t_4]. \eqno (3.30)
$$
By (3.28) $y_1(t)$ is a solution of the linear equation
$$
Y'(t) + \frac{y_1(t)}{p(t)} Y(t) + \sum\limits_{k=1}^n q_k(t) \exp\biggl\{\il{t}{\alpha_k(t)} \frac{y_1(\tau)}{p(\tau)} d \tau\biggr\} = 0, \phh t\in [t_1;t_4].
$$
Then by Cauchy formula
$$
y_1(t) = \exp\biggl\{ - \il{\xi}{t}\frac{y_1(\tau)}{p(\tau)} d\tau\biggr\}\biggl[y_1(\xi) - \il{\xi}{t}\Bigl(\sum\limits_{k=1}^n q_k(\tau) \exp\biggl\{\il{\xi}{\alpha_k(\tau)}\frac{y_1(s)}{p(s)} d s\biggr\}\Bigr)d\tau\biggr], \phh t\in [\xi;t_2],
$$
$t_1 \le \xi \le t_2$.
From here it follows that  $y_1(t) \ge \il{t}{t_2}\bigl(\sum\limits_{k\in \omega_+} q_k(\tau)\bigr) d\tau, \ph t\in [t_1;t_2].$
From here and from (3.30) we will get:
$$
\sum\limits_{k=1}^n q_k(t) \exp\biggl\{\il{t}{\alpha_k(t)} \frac{y_1(\tau)}{p(\tau)} d \tau\biggr\}\ge \sum\limits_{k\in\omega_+}^n q_k(t) \exp\biggl\{\il{t}{\alpha_k(t)} \frac{d\tau}{p(\tau)}\il{\tau}{t_2}\bigl(\sum\limits_{j\in\omega_+}q_j(s)\bigr)d s\biggr\} + \phantom{aaaaaaaaaaaaa}
$$
$$
\phantom{aaaaaaaaaaaaaaaaaaaaaaaa}+ \sum\limits_{k\in\omega_-}q_k(t)\frac{\il{t_1}{\alpha_k(t)}\frac{d\tau}{p(\tau)} + \varepsilon}{\il{t_1}{t}\frac{d\tau}{p(\tau)} + \varepsilon} \stackrel{def}{=} Q_\varepsilon(t), \phh t\in [t_1;t_2]. \eqno (3.31)
$$
Consider the Riccati equation
$$
y'(t) + \frac{1}{p(t)} y^2(t) + Q_\varepsilon(t) = 0, \phh t\in [t_1;t_2].
$$
Let $Y_\varepsilon(t)$ be the solution of this equation with  $Y_\varepsilon(t_1) = y_1(t_1)$.
Then by virtue of Lemma~3.1. from (3.31) it follows that  $Y_\varepsilon(t)$ exists on $[t_1;t_2)$  and
$$
Y_\varepsilon(t) \ge y_1(t), \phh t\in [t_1;t_2). \eqno (3.32)
$$
Then  $\phi_\varepsilon(t) \equiv \exp\biggl\{\il{t_1}{t}\frac{Y_\varepsilon(\tau)}{p(\tau)}d\tau\biggr\}, \ph t\in [t_1;t_2),$  is a solution of Eq. (2.9) on  $[t_1;t_2)$, which is continuable on  $[t_1;t_2]$
as a solution to Eq. (2.9).  From (3.32) it follows that  $\phi_\varepsilon(t)\ne~0, \linebreak t\in [t_1;t_2]$.
Hence (2.9) is not oscillatory on  $[t_1;t_2]$, which contradicts the condition c) of the theorem. Therefore it remains to suppose that the case II) takes place.  By virtue of   a) for this case the inequality
$$
\sum\limits_{k=1}^n q_k(t)\exp\biggl\{\il{t}{\alpha_k(t)}\frac{y_1(\tau)d\tau}{p(\tau)}\biggr\} \ge \sum\limits_{k\in \omega_2^-} q_k(t). \eqno (3.33)
$$
is fulfilled   on $[t_3;t_4]$.  Consider the Riccati equation
$$
y'(t) + \frac{1}{p(t)}y^2(t) + \sum\limits_{k\in \omega_2^-} q_k(t) = 0, \phh t\in [t_3;t_4].
$$
Let  $y_2(t)$  be the solution of this equation with $y_2(t_3) = y_1(t_3)$.
Then by virtue of Lemma~3.1 from  (3.33) it follows that $y_2(t)$ exists on  $[t_3;t_4)$  and
$$
y_2(t) \ge y_1(t), \phh t\in [t_3;t_4). \eqno (3.34)
$$
Therefore  $\phi_2(t) \equiv \exp\biggl\{\il{t_3}{t} \frac{y_2(\tau)}{p(\tau)} d\tau\biggr\}, \ph t\in [t_3;t_4)$, is a solution of Eq. (2.12), which can be continued on $[t_3;t_4]$,
as a solution to Eq. (2.12).  From (3.34) it follows that $\phi_2(t) \ne 0, \linebreak t\in [t_3;t_4]$. Therefore Eq. (2.12) is not oscillatory which contradicts the condition  d) of the theorem. The obtained contradiction completes the proof of the theorem.

\vskip 20pt

\centerline{\bf References}

\vskip 20 pt

\noindent
1. L. Berezansky, E. Braverman. Some Oscillation Problems For a Second Order Linear\linebreak \phantom{a} Delay  Differential Equations. Mathematical analysis and applications, 220,\linebreak \phantom{aa}pp. 719 - 740,  1998

\noindent
2. L. Berezansky and E. Braverman, Oscillation of a Second-Order Delay
Differential\linebreak \phantom{aa}Equation with Middle Term, Applied Mathematics Letters 13 (2000) 21 - 25.

\noindent
3. J. Dzurina, Oscillation Theorems for Second - Order Advanced Neutral Differential \linebreak \phantom{aa}Equations. Tatra Mt. Math. Publ. 48 (2011), 61 - 71.

\noindent
4. J. Ohriska, Oscillation of Second Order Linear Delay Differential Equations, Cont. Eur. \linebreak \phantom{aa}J. Math. 6 (3), 2008, 439 - 452.

\noindent
5. Z. Oplustil, J. Sramer, Some Oscillation Criteria for the Second - Order Linear Delay \linebreak \phantom{aa}Differential Equation. Mthematica Bohemica, Vol. 136 (2011), No 2, 195 - 204.

\noindent
6. T. Li, Yu. V. Rogovchenko and Ch. Zeng, Oscillation of Second - Order Neutral \linebreak \phantom{aa}Differential Equations,
Funccialaj Ekvacioj, 56 (2013), 111 - 120.

\noindent
7. W. B. Fite, Concerning the zeroes of the solutions of certain differential equations, \linebreak \phantom{aa}Trans. Amer. Math. Soc. 19 (1918), 341 - 352.

\noindent
8. A. I. Egorov, Uravneniya Rikkati (Riccatti Equations), Moscow, 2001.

\noindent
9. Ph. Hartman, Ordinary Differential Equations, SIAM, Classics in Applied Mathematics \linebreak \phantom{aa}38, Philadelphia 2002.

\end{document}